\documentclass[11pt]{article}
\usepackage{amsmath, amssymb, amsthm}

\theoremstyle{plain}

\newtheorem{thm}{Theorem}
\newtheorem{prop}[thm]{Proposition}
\newtheorem{lemma}[thm]{Lemma}

\theoremstyle{definition}
\newtheorem{defn}[thm]{Definition}
\newtheorem{remark}[thm]{Remark}

\newtheorem*{ackn}{Acknowledgement}

\begin{document}
\renewcommand{\today}{May 15, 2000}

\def\<{\langle} \def\>{\rangle}
\def\what{\widehat}
\def\Z{{\mathbb Z}}\def\N{{\mathbb N}} \def\C{{\mathbb C}}
\def\Q{{\mathbb Q}}\def\R{{\mathbb R}} \def\G{{\mathbb G}}
\def\rev{\mathop{\rm rev}\nolimits}
\def\Isom{\mathop{\rm Isom}\nolimits}

\def\Proof{\paragraph{Proof.}}
\def\Endproof{\endproof\medskip}
\def\notation{\paragraph{Notation.}}
\def\ackn{\paragraph{Acknowledgement.}}

\def\Mat{M} \def\inv{^{-1}} \let\bk\backslash
\def\Cal#1{{\cal#1}} \let\iff\Leftrightarrow
\let\liff\Longleftrightarrow \let\imply\Rightarrow

%
\def\al{\alpha}                 \def\be{\beta}
\def\ga{\gamma}
\def\ep{\epsilon}               \def\varep{\varepsilon}
%
%

\def\id{1}

\title{On the Toeplitz algebras of right-angled \\ and finite-type Artin
groups}

\author{
\hbox{John Crisp\footnote{Supported, during the preparation of this paper,
by an EPSRC
Research Assistantship at the University of Southampton, United Kingdom, and a Postdoctoral Grant from the Conseil Regional de Bourgogne, France.
The first
author also wishes to thank the Department of Mathematics at the University of
Newcastle, NSW, for their hospitality on numerous  occasions.}}
\hbox{\hskip.5cm and \hskip.5cm}
\hbox{Marcelo Laca\footnote{Supported by the Australian Research Council}}
}

\maketitle

\begin{abstract} The graph product of a family of groups lies somewhere
between their
direct and free products, with the graph determining which pairs of groups
commute and
which do not. We show that the graph product of quasi-lattice ordered groups is
quasi-lattice ordered, and, when the underlying groups are amenable, that it
satisfies  Nica's amenability condition for quasi-lattice orders
\cite{Nica}. As a consequence, the Toeplitz algebras of these groups are
universal for
covariant isometric representations on Hilbert space, and their
representations are
faithful if the isometries satisfy a properness condition given by Laca and
Raeburn
\cite{LR}.  An application of this to right-angled Artin groups gives a
uniqueness
theorem for the $C^*$-algebra generated by a collection of isometries such
that any
two of them either
$*$-commute or else have orthogonal ranges.  In contrast, the nonabelian
Artin groups
of finite type considered by Brieskorn and Saito \cite{BS} and Deligne
\cite{Del} have
canonical quasi-lattice orders that are not amenable in the sense of Nica,
so their
Toeplitz algebras are not universal  and the $C^*$-algebra generated by a
collection
of isometries satisfying the Artin relations fails to be unique.
\end{abstract}



\section{Introduction}\label{Intro}

Several celebrated results in $C^*$-algebra theory  assert that the
$C^*$-algebra generated by a semigroup of isometries does not depend on the
specific
isometries, provided they satisfy a properness condition. The situations
described by
these results are of considerable interest, stemming from the fact that the
algebraic
structure given by the semigroup operation determines a unique $C^*$-norm
on the
$*$-algebra generated by the isometries. As examples we have Coburn's
theorem on the
$C^*$-algebra generated by a single
 isometry \cite{cob}, Douglas's theorem on the $C^*$-algebra of a one parameter
semigroup of isometries \cite{dou}, and the generalization by Murphy to the
Toeplitz
$C^*$-algebra of a totally ordered group \cite{mur}; in all these cases the
properness
condition simply says that the isometries are not unitary.

Moving away from total orders on abelian groups, Nica \cite{Nica}
considered a class
of partially ordered groups
$(G,P)$ he called {\em quasi-lattice ordered}. Inspired by what happens
with the left
regular (Toeplitz) representation of the positive cone $P$, he isolated a key
covariance condition, which is automatic for total orders, and defined a
universal
$C^*$-algebra $C^*(G,P)$ whose representations are given by the covariant
isometric
representations of $P$. He proved that the uniqueness of the $C^*$-algebra
generated
by a covariant isometric representation depends on an amenability property
of the
quasi-lattice order that is strictly weaker than amenability of the
underlying group.
Indeed, he showed that Cuntz's result \cite{cun2} on the uniqueness of the
$C^*$-algebra $\mathcal T \mathcal O_n$ generated by $n$ isometries with
orthogonal
ranges can be seen as an amenability result for the canonical quasi-lattice
order on
the free group on $n$ generators. In this case the covariance condition
requires that
the generating isometries
 have orthogonal ranges, and the properness condition says that the sum of
these
ranges is not the whole Hilbert space.

In \cite{LR} Laca and Raeburn associated a semigroup dynamical system to each
quasi-lattice order and showed that the corresponding crossed product is
canonically
isomorphic to the universal $C^*$-algebra $C^*(G,P)$. This approach led to
two main
advances. The first one was the generalisation to all quasi-lattice orders
of some key
estimates of Cuntz
\cite{cun}, which provides a convenient framework in which to study
faithfulness of
representations and uniqueness properties.
 The second one was a direct proof of the amenability of the quasi-lattice
orders on a
large class of (nonamenable) free product groups, which widened the range of
application of the uniqueness results.

Direct products and free products of groups are both special cases of the
more general
construction of a graph product of groups (see Section
\ref{Sect2} below), and in this paper we address the  natural questions of
whether
graph products support quasi-lattice orders, and under which conditions
these are
amenable.
 Our main technical results are Theorem
\ref{graph-prods-qlo}, which shows that graph products of quasi-lattice
ordered groups
are indeed quasi-lattice ordered, and Theorem \ref{main-amen}, which gives a
sufficient condition for their amenability. Combining this with the results of
\cite{LR} in  Theorem
\ref{unique-gr-pr}, we characterize
 faithful representations and give a uniqueness result for the Toeplitz
algebras of
graph products.

An interesting class of examples is that of graph products of copies of
$(\Z,\N)$, otherwise known as graph groups or right-angled Artin groups. It
follows
from our main results that they are quasi-lattice ordered and amenable in
the sense of
\cite{Nica}, giving a unified statement of the amenability of the canonical
quasi-lattice orders on all free groups and all free abelian groups, as well as
providing many new examples of amenable quasi-lattice orders. The corresponding
Toeplitz $C^*$-algebras are thus universal and unique.
 We state this main result in terms of generators and relations in
Theorem~\ref{ra-artin}, which contains, as extreme cases,  Cuntz's theorem
(in which
the generating isometries have mutually orthogonal ranges) and a
multivariable version
of Coburn's theorem, (in which the generating isometries $*$-commute, i.e. they
commute with each other and with each other's adjoints). See \cite{sal} for
results
related to this last situation.

Other interesting quasi-lattice orders are provided by the family of Artin
groups of
finite type, with the embedded Artin monoid as positive cone \cite{BS,Del}.
These
examples, which include the braid groups, are {\em lattice groups}, because
every pair
of elements has a least common upper bound. In Section \ref{nonamen} we
prove, using
an argument essentially due to Nica, that if a group is lattice ordered and
amenable
as a quasi-lattice order, then the group itself has to be amenable. Thus,
in contrast
to what happens with graph products, only the Artin groups of finite type
that are
amenable (and hence abelian)
 give rise to amenable quasi-lattice orders. The nonabelian Artin groups of
finite
type appear then as an important class of groups having canonical
quasi-lattice orders
that are not amenable in the sense of \cite{Nica}. As a consequence, the
$C^*$-algebra
generated by a covariant isometric representation depends, in general, on
the specific
representation, Theorem~\ref{nonunique}.

\section{Graph products of groups}\label{Sect2}

Graph products were defined in the thesis of E.R.~Green \cite{Green}, and
have been
subsequently studied by various other authors.  We refer the reader to
\cite{HM} and
the references therein for further background.

Let $\Gamma$ denote a simplicial graph with vertex set $\Lambda$, and edge set
$E(\Gamma)\subseteq\{(I,J):I,J\in \Lambda \text{ and } I\neq J\}$. That is,
we assume
that $\Gamma$ has no edges which start and finish at the same vertex. We
say that
vertices $I$ and $J$ are \emph{adjacent} in $\Gamma$ if $(I,J)\in
E(\Gamma)$. Note
that a vertex is never adjacent to itself.  Given a family
$\{G_I\}_{I\in\Lambda}$ of
groups, we define the \emph{graph product} $\Gamma_{I\in\Lambda} G_I$ to be the
quotient of the free product
$*_\Lambda G_I$ by the smallest normal subgroup containing the elements
$x_1x_2x_1\inv x_2\inv$ for all pairs $x_1\in G_I$, $x_2\in G_J$ where
$I$ and $J$ are adjacent in $\Gamma$. When the $G_I$ are all copies of
$\mathbb Z$,
the graph product is called a {\em graph group}, or {\em right-angled Artin
group}. We
shall not need to assume that $\Gamma$ is finite.

Suppose, for the rest of this section, that we are given a graph $\Gamma$,
as above,
and groups $\{ G_I\}_{I\in\Lambda}$, and let $G=\Gamma_{I\in\Lambda} G_I$
denote the
graph product. We may take as a generating set for $G$ the set
\[\G =\coprod\limits_{I\in\Lambda}G_I\setminus\{ 1\}\,.\] Given $x\in\G$ we
write
$I(x)$ for the unique vertex $I$ such that $x\in G_I$. We say that $x$
\emph{belongs}
to $I(x)$.

By an \emph{expression} $X$ for an element $x\in G$ we mean a word in the
generators
$\G$ representing $x$. Given an expression $X=x_1x_2\cdots x_\ell$, the
elements $x_i$
are called \emph{syllables} of $X$ and the number $\ell$ is called the
\emph{length}
of $X$, written $\ell=\ell(X)$. We say that $I\in\Lambda$ is a \emph{vertex
of} $X$ if
$I=I(x_i)$ for $x_i$ a syllable of $X$.

Given an expression $X=x_1x_2\cdots x_\ell$ for $x\in G$, the graph product
relations
allow one to modify $X$ to obtain a different expression for
$x$ by replacing a subexpression $x_ix_{i+1}$ with $x_{i+1}x_i$ if
$I(x_i)$ is adjacent to $I(x_{i+1})$. In the terminology of \cite{HM}, such a
substitution is called a \emph{shuffle}, and we shall say that two
expressions are
\emph{shuffle equivalent} if one may be obtained from the other via a
finite sequence
of shuffles.  If the expression $X$ contains a subexpression of the form
$x_ix_{i+1}$,
with
$I(x_i)=I(x_{i+1})$, then we may give a shorter expression for $x$ by an
\emph{amalgamation}, that is by deleting  $x_ix_{i+1}$ in the case that
$x_{i+1}=x_i\inv$ or otherwise by replacing the two syllables $x_ix_{i+1}$
with the
single syllable $\what x_i\in\G$ such that $\what x_i=x_ix_{i+1}$.

We say that an expression is \emph{reduced} if it is not shuffle equivalent
to any
expression which admits an amalgamation.

\begin{lemma}\label{Xred} Given an expression $X=x_1x_2...x_\ell$, the
following are
equivalent:
\begin{description}
\item{(i)} $X$ is reduced;
\item{(ii)} for all $i<j$ such that $I(x_i)=I(x_j)$ there exists $k$ such
that $i<k<j$
and $I(x_k)$ is not adjacent to $I(x_i)$.
\end{description}
\end{lemma}

\Proof That (ii) is a consequence of (i) is obvious. We note that the truth of
condition (ii) is invariant under shuffles, while an amalgamation can never be
performed on an expression satisfying (ii). Thus (ii)$\imply$(i).\Endproof

Given an expression $X=x_1x_2..x_\ell$ for $x\in G$, one may easily produce
a reduced
expression for $x$ via a finite process of shuffling and amalgamating. If
$X$ fails
condition (ii) of Lemma \ref{Xred}, that is there exist $i<j$ such that
$I(x_i)=I(x_j)$ and $I(x_i)$ commutes with $I(x_k)$ for all $i<k<j$, then
an obvious
sequence of shuffles will allow the syllables $x_i$ and $x_j$ to be
amalgamated,
reducing the length of $X$. One may continue reducing the length in this
manner until
condition (ii) is satisfied. The following theorem, due to Green \cite{Green},
reduces the solution of the word problem in a graph product of groups to
the solution
of the word problem in each of the component groups.(see \cite{HM} for an
equivalent
formulation).

\begin{thm}[E.R.Green \cite{Green}]\label{NormalForm} Any two reduced
expressions for
the same element of $G$ are shuffle equivalent.  \end{thm}

As a simple consequence of Theorem \ref{NormalForm} we may define the
\emph{length}
$\ell(x)$ of an element $x\in G$ to be the length of any reduced expression
$X$ for
$x$, note that $\ell(x)$ is also the minimal length of any expression for
$x$. If $X$
is a reduced expression for $x$, we are also justified in referring to a
syllable or a
vertex of $X$ as a \emph{syllable of $x$} or a \emph{vertex of $x$},
respectively.

Let $X=x_1x_2...x_\ell$ be a reduced expression, and $x_i$ any syllable of
$X$. If
$I(x_i)$ is adjacent to $I(x_j)$ for all $j<i$ then we say that $x_i$ is an
\emph{initial syllable} of $X$, and that $I(x_i)$ is an \emph{initial
vertex} of $X$.
We write $\Delta(X)$ for the set of all initial vertices of $X$. The
following facts
are easily checked:

\begin{lemma}\label{lemm1} Let $X=x_1x_2...x_\ell$ be a reduced expression.
\begin{description}
\item{(i)} If $x_i$ is an initial syllable of $X$ then $X$ is shuffle
equivalent to
the expression \[ x_ix_1..x_{i-1}x_{i+1}..x_\ell\,.\]
\item{(ii)} The initial vertices of $X$ are pairwise adjacent.
\item{(iii)} For each $I\in\Delta(X)$ there exists a unique initial
syllable of $X$
belonging to $I$. We define the function $\delta_X:\Delta(X)\to\G$ such that
$\delta_X(I)$ is equal to the initial syllable of $X$ belonging to $I$.
\item{(iv)} If $X'$ is another expression which is shuffle equivalent to
$X$ then
$\Delta(X')=\Delta(X)$ and $\delta_{X'}=\delta_X$.
\end{description}
\end{lemma}

By virtue of Lemma \ref{lemm1}(iv) and Theorem \ref{NormalForm} we may make the
following definitions.

\begin{defn}\label{xI's} Let $x\in G$. We define the set of \emph{initial
vertices} of
$x$ to be the set $\Delta(x)=\Delta(X)$ where $X$ is any reduced expression
for $x$.
We also define, for each $I\in\Lambda$, an element $x_I\in G_I$ as follows.
Take any
reduced expression $X$ for $x$, and set
\[ x_I =
\begin{cases}
\delta_X(I)\qquad&\text{if }I\in\Delta(x)\,,\\ 1\qquad&\text{if
}I\notin\Delta(x)\,.
\end{cases}
\] Given an expression $X=x_1x_2..x_\ell$ let $\rev(X)$ denote the
expression $x_\ell
x_{\ell-1}...x_1$. This allows us to define the set of \emph{final
vertices} of $x\in
G$ to be the set $\Delta^r(x)=\Delta(\rev(X))$ where $X$ is any reduced
expression for
$x$. Similarly also, for each $I\in\Lambda$, we define the element
\[ x^r_I =
\begin{cases}
\delta_{\rev(X)}(I)\qquad&\text{if }I\in\Delta(x)\,,\\ 1\qquad&\text{if
}I\notin\Delta(x)\,,
\end{cases}
\] where $X$ is any reduced expression for $x$.
\end{defn}

Observe that, by using Lemma \ref{lemm1}(i) and (ii), one may always find,
for any
given element $x\in G$, a reduced expression  which begins with the
product, in any
order, of the $x_I$ for $I\in\Delta(x)$. Equally, one may always find a reduced
expression for $x$ which ends with the product, in any order, of the
$x^r_I$ for
$I\in\Delta^r(x)$.

\begin{lemma}\label{XZY} Given $x,y\in G$, let
$D=\Delta^r(x)\cap\Delta(y)$, and
suppose that $z_I=x^r_Iy_I$ is nontrivial for each $I\in D$. We define an
expression
$Z=\prod_{I\in D}z_I$ (in any order). Then, if $X \cdot (\prod_{I\in D}
x^r_I)$ and
$(\prod_{I\in D} y_I) \cdot Y$ are reduced expressions for $x$ and $y$
respectively,
we have that $XZY$ is a reduced expression for $xy$.
\end{lemma}

\Proof Since each element $z_I$ is nontrivial, it is clear that both
expressions $XZ$
and $ZY$ are reduced since they are formally equivalent to the given reduced
expressions for $x$ and $y$ respectively.  (Two expressions
$x_1x_2..x_\ell$ and
$y_1y_2..y_k$ are \emph{formally equivalent} if $k=\ell$ and
$I(x_i)=I(y_i)$ for all
$1\leq i\leq \ell$.)

Suppose that $XZY=w_1w_2..w_\ell$ is not reduced.  Then, by Lemma
\ref{Xred}, one may
find $i<j$ such that $I(w_i)=I(w_j)$ and $I(w_i)$ is adjacent to $I(w_k)$
for all $k$
with $i<k<j$. But $w_i$ and $w_j$ are not both in $XZ$ nor both in $ZY$,
since these
are reduced expressions. Thus $w_i$ must be from the subexpression $X$ and
$w_j$ from
$Y$. But now it follows that $w_i$ is a final syllable of $XZ$ and $w_j$ an
initial
syllable of $ZY$.  That is $I(w_i)=I(w_j)=J$ for some $J\in D$. But this
contradicts
the fact that there is a syllable $w_k=z_J$ lying between $w_i$ and $w_j$.
\Endproof

\section{Quasi-lattice orders and their graph products}\label{Sect3}

Let $G$ be a group and let $P$ be a submonoid (subsemigroup containing the
identity)
of $G$ such that $P\cap P\inv=\{\id\}$.  Then we may define a
left-invariant partial
order on $G$ by $x\leq y$ whenever
$x\inv y\in P$.  Note that $x\in P$ if and only if $1\leq x$.  We observe,
indeed,
that every left-invariant partial order on $G$ arises in this fashion.  We
say that
$(G,P)$ is a \emph{partially ordered group} with \emph{positive cone} $P$.

\begin{defn}\label{defqlo} A partially ordered group $(G,P)$ is
\emph{quasi-lattice ordered} if every finite set in $G$ with an upper bound
in $G$ has
a (necessarily unique) least upper bound in $G$. Equivalently, every pair
$x$, $y$ of
elements of $G$ with a common upper bound in $G$ has a least upper bound,
which we
denote by $x\vee y$.  If $x$ and
$y$ have no common upper bound in $G$, then we write $x\vee y=\infty$ for
convenience.
 \end{defn}

Given the group $G$ with positive cone $P$ we may equally define a
right-invariant
partial order on $G$ by $x\leq_r y$ whenever $yx\inv\in P$. If $x$ and $y$
have a
greatest lower bound for $\leq_r$, we denote it by
$x\wedge_r y$. Clearly one has
\begin{equation}\label{eqn1} x\leq_r y \ \text{ if and only if }\
y\inv\leq x\inv,
\quad  \text{ and }\quad x \wedge_r y = ( x\inv \vee y\inv )\inv.
\end{equation}

\begin{lemma}\label{qlo} For a partially ordered group $(G,P)$ the following
statements are equivalent.

\begin{enumerate}

\item[(i)] $(G,P)$ is a quasi-lattice order.

\item[(ii)] Every finite set in $G$ with a common upper bound in $P$ has a
least upper
bound in $P$.

\item[(iii)] Every element $x$ of $G$ having an upper bound in $P$ has a
least upper
bound in $P$.

\item[(iv)] If $x \in P P\inv$ then there exist a pair of elements $a,b\in
P$ with
$x=ab\inv$ and such that for every
$u,v\in P$ with $a b\inv = uv\inv$, one has $a\leq u$ and $b\leq v$. (The
pair $a,b$
is clearly unique.)

\item[(v)] Every pair $u,v$ of elements in $P$ has a greatest lower bound
$u\wedge_r
v$ with respect to the right-invariant partial order on
$G$.

\item[(vi)]  If $x \in P P\inv$ then there exist a pair of elements $a,b\in
P$ with
$x=ab\inv$ and such that
$a\wedge_r b = 1$.
\end{enumerate}
\noindent Assuming that (i)--(vi) hold, and given $x\in PP\inv$, there is
in fact a
unique pair $a,b\in P$ satisfying statement (vi), being precisely the pair
$a,b$ of
statement $(iv)$.
\end{lemma}

\Proof We prove (i)$\implies$(ii)$\implies$(iii)$\implies$(i) first.
Clearly (ii)
follows from (i) by observing that a (least) upper bound in $P$ for a finite set $F$ is the same thing as a (least) upper bound for $F\cup\{ 1\}$. Statement
(iii) is obviously a special case of (ii).  Finally, (iii) states
that the
condition of Definition \ref{defqlo} holds for all pairs $x,1$, and (i) may be
recovered from this by left invariance of the partial order.  \smallskip

(iii)$\implies$(iv): Suppose $x \in P P\inv$.  Then
$x$ has an upper bound in $P$, and by (iii) we can take $a$ to be the least
upper
bound in $P$, and write $x=ab\inv$ for $b\in P$. If $x=uv\inv$ for
$u,v\in P$, then
$u$ is an upper bound for $x$ so we have $a\leq u$, and hence also $b\leq
v$ (since
$b\inv v = a\inv u\in P$).  \smallskip

(iv)$\implies$(v): Given $u,v\in P$, let $x=uv^{-1}$. Then the pair
$a,b\in P$ of (iv) determines an element $e:=a\inv u=b\inv v$ which one
checks to be
the greatest right lower bound $u\wedge_r v$.  For if $w$ is a common right
lower
bound then $u=cw$ and $v=dw$, for some $c,d\in P$, and one has
$cd\inv=uv\inv =ab\inv
$.  But, by (iv), both $a\leq c$ and $b\leq d$, so that $w\leq_r e$.
\smallskip

(v)$\implies$(vi): Given $x\in PP^{-1}$, choose $u,v \in P$ such that
$x=uv^{-1}$.
Using (v) we may write $ u = a (u\wedge_r v)$ and $v = b (u\wedge_r v)$ for
some $a,b
\in P$. Clearly $uv\inv = ab\inv $ and $a \wedge_r b = 1$. \smallskip

(vi)$\implies$(iii): Let $x \in G$ have an upper bound $u$ in $P$.  Then $x
\in P
P\inv$, and by (vi) we can write $ x = a b\inv$ for $a,b \in P$ with
$a\wedge_r b
=1$.  Since $a b\inv \leq u$, we have that $b a\inv u \in P$, in other
words, $u\inv a
\leq_r b$.  It is obvious that $u\inv a \leq_r a$, so we must have that
$u\inv a \leq_r
a\wedge_r b = 1$.  Hence $a\inv u \in P$ and $a \leq u$.  So we have that
$a$ is the
least upper bound of $x$ in $P$.
\smallskip

Finally suppose that (i)--(vi) hold and let $x=uv\inv$ for $u,v\in P$.
Then, as in
the proof of (iv)$\implies$(v), we have $u\wedge_r v=a\inv u=b\inv v$ where
$a,b\in P$
is the unique pair of statement (iv).  Thus the pair $u,v$ satisfies the
condition of
statement (vi), namely
$u\vee_r v=1$, if and only if $u=a$ and $v=b$.  \Endproof

\begin{remark}
 We make the following remarks concerning Lemma \ref{qlo}.

\begin{itemize}
\item[(1)] Our definition of quasi-lattice order differs slightly from the one
originally given in \cite[\S 2.1]{Nica}, which appears here as statement
(ii). Nica
had also given an equivalent form of (ii) consisting of two conditions:
QL1, which is
statement (iii), and QL2, which is the statement of Definition \ref{defqlo}
for all
pairs $x,y\in P$. By Lemma \ref{qlo} the various definitions are equivalent; in
particular QL2 is not needed as it follows from QL1.

\item[(2)] In \cite{Nica} and \cite{LR} least upper bounds are always
assumed to be in
$P$, but no such assumption is made here. This causes a slight discrepancy
in notation:
in this paper, the least upper bound in $P$ of $x$ and $y$ would be written
$x\vee
y\vee 1$.

\item[(3)] While statements (v) and (vi) may appear to be conditions only
on the
monoid $P$, in fact they are not, because one must have that
$ w \leq_r u\wedge_r v$ not only for every common right lower bound $w$ in
$P$ but
also for every common right lower bound $w$ in $G$.

\end{itemize} \end{remark}

Suppose, now, that $G=\Gamma_{I\in\Lambda}G_I$ is a graph product in which
each group
$G_I$, for $I\in\Lambda$, is partially ordered with positive cone $P_I$.
We say that
a reduced expression $X=x_1x_2...x_\ell$ is \emph{positive} if $x_i\in
P_{I(x_i)}$ for
all $i=1,2,...,\ell$. Note that this property is invariant under shuffle
equivalence.
We say  that an element $x\in G$ is \emph{positive} if it has a positive
reduced
expression.  It follows, by Theorem \ref{NormalForm}, that every reduced
expression
for a positive element is positive. Let $P$ denote the submonoid of $G$
consisting of
all positive elements. This is just the submonoid generated by the union of
the $P_I$
for $I\in \Lambda$.  
Moreover, two positive elements are equal in $P$ if and only if their
reduced expressions are shuffle equivalent. Thus $P$ may be presented as the 
\emph{monoid graph product} $\Gamma_{I\in\Lambda} P_I$, that is the monoid 
obtained from the free product $*_{I\in\Lambda}P_I$ by introducing the 
relations $xy=yx$ for all $x\in P_I$ and $y\in P_J$ with $I$ adjacent to $J$ 
in $\Lambda$.

It is easily seen that $P\cap P\inv =\{1\}$ and hence that
$(G,P)$ is a partially ordered group.  We refer to $(G,P)$  as the
\emph{graph product
over $\Gamma$ of the partially ordered groups} $\{
(G_I,P_I)\}_{I\in\Lambda}$, and
write $(G,P)=\Gamma_{I\in\Lambda}(G_I,P_I)$.  Note that each $(G_I,P_I)$ is a
partially ordered subgroup of $(G,P)$.  That is, the inclusion map is order
preserving.
\medskip

In what follows, we suppose that $(G,P)=\Gamma_{I\in\Lambda}(G_I,P_I)$ is a
graph
product of quasi-lattice ordered groups.  We will use the notation
$\Delta^r(x,y)$ for
the intersection $\Delta^r(x)\cap\Delta^r(y)$ of the final vertex sets of
elements $x$
and $y$ in $G$.

\begin{lemma}\label{red-quot} Given any pair $u,v\in P$ there exist $a,b\in P$
satisfying the following conditions:
\begin{enumerate}
\item[(i)] { $ab\inv =uv\inv$, with
$a\leq u$, $b\leq v$, and $a_I^r\wedge_r b_I^r=1$ for all $I\in\Delta^r(a,b)$.

\noindent (Note that by Lemma \ref{qlo} (v) $\wedge_r$ is defined in each
quasi-lattice order $(G_I,P_I)$.)  }

\item[(ii)]{ Writing reduced expressions $A\cdot(\prod_{\Delta^r(a,b)}
a^r_I)$ and
$(\prod_{\Delta^r(a,b)} (b^r_I)\inv)\cdot B$ for $a$ and
$b\inv$ respectively, one has that
\begin{equation}\label{red-expr} A\cdot
\bigl\lbrack\prod_{\Delta^r(a,b)}(a_I^r(b^r_I)\inv)\bigr\rbrack\cdot B
\end{equation} is a reduced expression for $ab\inv$.}
\end{enumerate}
\end{lemma}

\Proof (i): We proceed by induction on $\ell(u)+\ell(v)$, the case where
$u=v=1$ being trivially true. Given $u,v\in P$ we have

\[ u=\what{u}\cdot\prod_{I\in\Delta^r(u,v)}(u_I^r\wedge_r v_I^r)
\hskip5mm\text{ and }\hskip5mm
v=\what{v}\cdot\prod_{I\in\Delta^r(u,v)}(u_I^r\wedge_r
v_I^r)\,,
\]

\noindent where $\what{u},\what{v}\in P$ and
$\what{u}\what{v}\inv=uv\inv$.  Now if both $u_I^r(u_I^r\wedge_r
v_I^r)\inv$ and
$v_I^r(u_I^r\wedge_r v_I^r)\inv$ are nontrivial for each
$I\in\Delta^r(a,b)$, then
these are precisely the final syllables of $\what{u}$ and $\what{v}$
respectively (in
particular,
$\Delta^r(\what{u},\what{v})=\Delta^r(u,v)$).  Putting $a=\what{u}, b=\what{v}$
satisfies the claim in this case.  Otherwise, some
$u_I^r\wedge_r v_I^r$ equals either $u_I^r$ or $v_I^r$ for some
$I\in\Delta^r(a,b)$.  In this case $\ell(\what{u})+\ell(\what{v})
<\ell(u)+\ell(v)$,
and the result follows by induction.  \smallskip

 (ii): This is a straightforward application of Lemma \ref{XZY}, the condition
$a_I^r\wedge_r b_I^r=1$ ensuring that the syllables
$(a_I^r(b_I^r)\inv)$ are nontrivial.\Endproof

\begin{thm}\label{graph-prods-qlo} A graph product
$(G,P)=\Gamma_{I\in\Lambda}(G_I,P_I)$ of quasi-lattice ordered groups is a
quasi-lattice ordered group.
\end{thm}

\Proof We prove that $(G,P)$ satisfies condition (iv) of Lemma \ref{qlo}.

Given $x\in PP\inv$, Lemma \ref{red-quot} implies that $x$ has a reduced
expression
$X=x_1x_2...x_m$ in which each syllable $x_i$ lies in
$P_{I(x_i)}P_{I(x_i)}\inv$.  (More specifically, $X$ may take the form of
Equation
(\ref{red-expr}), with the strictly positive syllables appearing first and
strictly
negative syllables appearing last).  Since each
$(G_I,P_I)$ is quasi-lattice ordered we may, in view of Lemma
\ref{qlo}(vi), write each $x_i$ uniquely as $a_ib_i\inv$ with
$a_i\wedge_r b_i=1$.  Define elements ${\bf a}$, ${\bf b}\in P$ by ${\bf
a}=a_1a_2...a_m$ and ${\bf b}=b_mb_{m-1}...b_1$. Since any shuffle of the
syllables
$x_1,x_2,..,x_m$ induces a legal shuffle of the nontrivial $a_i$'s (resp. the
nontrivial $b_i$'s) it follows by Theorem \ref{NormalForm} that the pair of
elements
${\bf a, b}$ is uniquely determined by $x$ (independently of the choice of
reduced
expression).

Given any pair $u,v\in P$ with $uv\inv=x$, let $a,b\in P$ be any pair as in
Lemma
\ref{red-quot}.  The unique pair ${\bf a,b}$ which we have just defined may be
computed from any reduced expression for
$x$, in particular from expression (\ref{red-expr}) of Lemma
\ref{red-quot}(ii).  It follows that ${\bf a}=a$ and ${\bf b}=b$.
  This shows that $x={\bf ab}\inv$ and that ${\bf a}\leq u$ and
${\bf b}\leq v$ for every $u,v\in P$ such that $x=uv\inv$, as required by Lemma
\ref{qlo}(iv).
\Endproof

\begin{lemma}\label{first} Suppose that $(G,P)$ is a graph product of partially
ordered groups $(G_I,P_I)$ for $I\in\Lambda$. Let $x,z\in P$ be such that
$1\leq x\leq
z$. Then, for each $I\in\Lambda$,
\begin{enumerate}
\item[(i)] $x_I\leq z_I$, \hskip5mm and
\item[(ii)] writing $x=x_Ix'$, either $x_I=z_I$ or $I$ is adjacent to every
vertex of
$x'$.
\end{enumerate}
\end{lemma}

\Proof  We have $z=xw$ where $x,z$ and $w$ are all positive. Let
$D=\Delta^r(x)\cap\Delta(w)$ and take reduced expressions $X\cdot
(\prod_{J\in D}
x^r_J)$ for $x$ and $(\prod_{J\in D} w_J)\cdot W$ for $w$. Note that
$x^r_Jw_J$ is
nontrivial, for each $J\in D$, since it is a product of nontrivial positive
elements.
Thus, by Lemma \ref{XZY}, we have the following (necessarily positive) reduced
expression for $z$ :
\[ Z = X\cdot (\prod\limits_{J\in D} (x^r_Jw_J))\cdot W\ .
\]

If $z_I=1$ then the Lemma is trivially true. Thus we suppose $z_I\neq 1$
and ask where
the initial syllable $z_I$ might appear in the above expression.
\medskip

\noindent\emph{Case 1.} If the syllable $z_I$ falls in the subexpression
$X$, then it
is also an initial syllable of the reduced expression given for $x$. Thus
$x_I=z_I$,
satisfying (i) and (ii).
\medskip

\noindent\emph{Case 2.} If $I\in D$ and the initial syllable $z_I$ happens
to be the
syllable $(x^r_Iw_I)$, then every vertex of $X$ and every $J\in D$ other
than $I$
itself is adjacent to $I$. It follows that $x^r_I$ is also an initial
syllable of $x$.
That is $x_I=x^r_I$, and $x'$ has a reduced expression
\[ X\cdot (\prod\limits_{J\in D;J\neq I} x^r_J)\,.
\]
So  $z_I=x_Iw_I$ giving part (i), and $I$ is adjacent to every vertex
of $x'$
giving part (ii).
\medskip

\noindent\emph{Case 3.} If the initial syllable $z_I$ falls in the
subexpression $W$,
then every vertex of $x$ is adjacent to $I$ and, in particular, $x_I=1$.
Both parts of
the Lemma are again satisfied.
\Endproof

\begin{defn}\label{I-adj} We consider a graph product
$(G,P)=\Gamma_{J\in\Lambda}
(G_J,P_J)$ of quasi-lattice orders, and choose $I\in\Lambda$. Given
elements $x,y\in
P$, write $x=x_Ix'$ and $y=y_Iy'$. We say that the elements $x,y\in P$ are
\emph{$I$-adjacent} if the following three conditions hold:
\begin{enumerate}
\item[(a)] $x_I\vee y_I\neq\infty$;
\item[(b)] either $x_I=x_I\vee y_I$ or $I$ is adjacent to every vertex of $x'$;
\item[(c)] either $y_I=x_I\vee y_I$ or $I$ is adjacent to every vertex of $y'$.
\end{enumerate}
\end{defn}

This definition allows us to give an inductive algorithm for deciding
whether two
elements $x,y\in P$ have a common upper bound, and for computing $x\vee y$
when it
exists:

\begin{prop}\label{lubs} Suppose that
$(G,P)=\Gamma_{J\in\Lambda}(G_J,P_J)$ is a graph product of quasi-lattice
ordered
groups. Let $x,y\in P$ and, for an arbitrary choice of $I\in\Lambda$, write
$x=x_Ix'$
and $y=y_Iy'$. Then we have the following.
\begin{enumerate}
\item[(i)] The elements $x,y\in P$ have a common upper bound if and only if
they are
$I$-adjacent and $x'\vee y'\neq\infty$.
\item[(ii)] Suppose that the elements $x,y\in P$ do have a common upper
bound. Then
\[ x\vee y = (x_I\vee y_I)\cdot(x'\vee y')\,.
\] Note that conditions (b) and (c) of Definition \ref{I-adj} apply to this
expression.
\end{enumerate}
\end{prop}

\Proof Suppose initially that the elements $x$ and $y$ are $I$-adjacent
and, by
condition (a) of Definition \ref{I-adj}, write $x_I\vee y_I=x_Iu=y_Iv$ for some
$u,v\in P_I$. Then, by condition (b), either $u=1$ or $I$ is adjacent to
every vertex
of $x'$. In either case one has that $ux'=x'u=u\vee x'$. Similarly, by
condition (c),
we have $vy'=y'v=v\vee y'$. Therefore, supposing in addition that $x'\vee
y'\neq\infty$ we have:
\begin{equation}\label{calculation}
\begin{aligned} (x_I\vee y_I)(x'\vee y')&= x_Iux'\vee y_Ivy'\\
				&= x_I(u\vee x')\vee y_I(v\vee y')\\
				&= (x_I\vee y_I)\vee x\vee y \\
				&= x\vee y
\end{aligned}
\end{equation}

\noindent It follows that $x$ and $y$ have a common upper bound and
$x\vee y$ is given by the above equation.

On th
e other hand, suppose that $z=x\vee y\neq\infty$. We first observe
that $x$ and
$y$ are $I$-adjacent. Clearly $x_I\vee y_I\leq z$, justifying (a), and
conditions (b)
and (c) follow directly from Lemma \ref{first}(ii). Finally, since $x\vee
y\neq\infty$, the equalities of Equation (\ref{calculation}) may be read in
reverse
order to show that $x'\vee y'\neq\infty$.
\Endproof

\begin{remark} To see that Proposition \ref{lubs} leads to an effective
algorithm, we
note that by always choosing $I\in\Delta(x)\cup\Delta(y)$ we ensure that
$\ell(x')+\ell(y')$ is always strictly less than $\ell(x)+\ell(y)$.
\end{remark}

\section{Amenability for graph products of quasi-lattice orders}

Recall from \cite{Nica} that an isometric representation $V: P \to
\Isom(\mathcal H)$ on a Hilbert space $\mathcal H$ is {\em covariant} if it is
compatible with the quasi-lattice structure in the sense that
\[ V_x V_x^* V_y V_y^* = V_{x\vee y} V_{x\vee y}^* \qquad \text{for } x,y
\in P.
\] The notation is meant to include the convention $V_\infty = 0$, so in
particular
covariance implies $V_x V_x^* V_y V_y^* = 0$ when $x$ and $y$ have no
common upper
bound.

The main example of such a representation is the Toeplitz representation
$T: P \to
\Isom(\ell^2(P))$, defined by $T_x
\varepsilon_y := \varepsilon_{xy}$, where $\varepsilon_x$ denotes the typical
orthonormal basis vector of $\ell^2(P)$.
 The $C^*$-algebra generated by the
$T_x$ is called the Toeplitz $C^*$-algebra of the quasi-lattice ordered
group $(G,P)$
and is denoted $\mathcal T(G,P)$. Nica also considered the  $C^*$-algebra
$C^*(G,P)$,
universal for covariant isometric representations of $P$ and made the following
definition.
\begin{defn} When the canonical homomorphism of
$C^*(G,P)$ to $\mathcal T(G,P)$ is injective we say that
$(G,P)$ is an {\em amenable quasi-lattice order}.
\end{defn}
There is a semigroup $C^*$-dynamical system $(B_P,P,\alpha)$ canonically
associated to
$P$, in which $B_P$ is the $C^*$-subalgebra of
$\ell^\infty(P)$ generated by the characteristic functions $1_y$ of the
semi-infinite
intervals $[y,\infty)$ for $y\in P$; the endomorphism
$\alpha_x$ corresponding to $x\in P$ is defined by $\alpha_x(1_y) = 1_{xy}$.
Covariant isometric representations of $P$ are in one to one correspondence
with
covariant representations of the semigroup dynamical system
$(B_P,P,\alpha)$ and this
leads to the realisation of $C^*(G,P)$ as a semigroup crossed product, see
\cite[Section 2]{LR} for the details.  There is a canonical conditional
expectation
from $B_P \rtimes_\alpha P$ onto $B_P$, which is faithful if and only if
$(G,P)$ is
amenable \cite[\S 4.3]{Nica}. This property, taken as the definition of
amenability of
$(G,P)$ in \cite{LR}, is instrumental in the direct proof of amenability
for free
product orders, which we aim to generalise in this section.

We will need the following reformulation of Proposition 6.6 of \cite{LR}.

\begin{prop}[Laca, Raeburn \cite{LR}]\label{LacaRae}  Suppose
$\phi:(G,P)\to (\Cal
G,\Cal P)$ is an order preserving homomorphism of quasi-lattice ordered
groups such
that, whenever $x,y \in P$ have a common upper bound in
$P$,
\begin{enumerate}
\item[(a)]\ \  $\phi(x)=\phi(y)$ only if $x=y$, and
\item[(b)]\ \  $ \phi(x)\vee\phi(y) = \phi(x\vee y)$.
\end{enumerate}
 If $\Cal G$ is an amenable group, then $(G,P)$ is an amenable
quasi-lattice order.
\end{prop}

\begin{remark} Proposition 6.6 of \cite{LR} should have been stated like this.
 The reason is that the proof indicated there, modelled on that of
\cite[Proposition 4.2]{LR}, requires that the conditional expectation for
the coaction
of $\Cal G$ on $C^*(G,P)$ be faithful, which is true if $\mathcal G$ is
amenable by
\cite{quigg}. We do not know whether Proposition 6.6 of \cite{LR} is correct as
originally stated; however all that is required for the other results in
\cite{LR} is
the version stated above.
\end{remark}

Suppose now that $(G,P)=\Gamma_{I\in\Lambda}(G_I,P_I)$ is a graph product of
quasi-lattice ordered groups. We define the group homomorphism
\[
\phi:(G,P)\longrightarrow \bigoplus\limits_{I\in\Lambda} (G_I,P_I)
\] such that each factor $(G_I,P_I)$ of $(G,P)$ is mapped to the
corresponding factor
in $\bigoplus_{I\in\Lambda} (G_I,P_I)$ via the identity on $G_I$. In what
follows we
shall, for simplicity, write $\phi(u)$ as $u$ whenever $u\in G_I$ for some
$I\in\Lambda$.

We view the direct product $\bigoplus_{I\in\Lambda} (G_I,P_I)$ as a graph
product
(over the full graph on $\Lambda$). Let $x\in G$ and let $X$ be any reduced
expression
for $x$. Then, choosing a vertex $I\in\Lambda$, we observe that $\phi(x)_I$
is simply
the product of all those syllables of $X$ which belong to the vertex $I$,
taken in the
order in which they appear.  In particular, if $x=x_Ix'$ then
$\phi(x)_I=x_I\phi(x')_I$. On the other hand, $\phi(x)_J=\phi(x')_J$ for
all $J\neq I$.

\begin{lemma}\label{lemA}
Let $\phi$ be the map defined above. Suppose that $x,y\in P$ satisfy $x\vee
y\neq\infty$ and, for an arbitrary choice of $I\in\Lambda$, write $x=x_Ix'$ and
$y=y_Iy'$. Then
\begin{equation}\label{lemAeqn}
\phi(x)_I\vee\phi(y)_I = (x_I\vee y_I)\cdot(\phi(x')_I\vee\phi(y')_I)\,.
\end{equation}
\end{lemma}

\Proof We refer to the final claim of Proposition \ref{lubs}, which states
that either
$y_I\leq x_I=x_I\vee y_I$ or $I$ is adjacent to every vertex of $x'$. The
latter
condition implies, in particular, that $\phi(x)_I=x_I$ and $\phi(x')_I=1$.
A similar
statement also holds with respect to $y'$. Thus we have the following four
cases to
consider:
\medskip

\noindent\emph{Case 1.} $\phi(x)_I=x_I$ and $\phi(y)_I=y_I$: Since in this case
$\phi(x')_I$ and $\phi(y')_I$ are both trivial, Equation (\ref{lemAeqn}) is
self-evident.
\medskip

\noindent{\emph{Case 2.} $\phi(x)_I=x_I$ and $x_I\leq y_I=x_I\vee y_I$}: In
this case
$\phi(x)_I\leq y_I\leq \phi(y)_I$. Also $\phi(x')_I=1$ and $x_I\vee
y_I=y_I$. Thus
Equation (\ref{lemAeqn}) reduces to $\phi(y)_I=y_I\phi(y')_I$.
\medskip

\noindent{\emph{Case 3.} $y_I\leq x_I=x_I\vee y_I$ and $\phi(y)_I=y_I$}:
This case is
similar to Case 2.
\medskip

\noindent{\emph{Case 4.} $x_I=y_I=x_I\vee y_I$}: Equation (\ref{lemAeqn})
follows in
this case by left invariance of the quasi-lattice order.
\Endproof

\begin{prop}\label{mainprop} Let $(G,P)=\Gamma_{I\in\Lambda}(G_I,P_I)$ be a
graph
product of quasi-lattice ordered groups. Then the map
$\phi:(G,P)\to\bigoplus_{I\in\Lambda}(G_I,P_I)$ (defined by the identity on
each
factor) is an order preserving homomorphism such that, whenever $x,y \in P$
have a
common upper bound in $P$, the following hold:
\begin{enumerate}
\item[(a)]\ \  $\phi(x)=\phi(y)$ only if $x=y$, and
\item[(b)]\ \  $ \phi(x)\vee\phi(y) = \phi(x\vee y)$.
\end{enumerate}
\end{prop}

\Proof The induced map is clearly an order preserving homomorphism. Suppose
throughout
that $x,y\in P$ have a common upper bound. \medskip

We first prove condition (a) assuming that (b) holds: Observe first of all
that if
$u\in P$ then  $u_I\leq \phi(u)_I$ for every $I$.
So  $\phi(u)=1$ implies that $u=1$.  Suppose that
$\phi(x)=\phi(y)$. Then, by condition (b), we have
\[
\phi(x\vee y)=\phi(x)\vee\phi(y)=\phi(x)=\phi(y)\,.
\] Writing $x\vee y=xu=yv$ for $u,v\in P$, it follows that
$\phi(u)=\phi(v)=1$. But by
the preceding observation we must then have $u=v=1$ and hence $x=y$.
\medskip

We now prove condition (b): Choose $I\in \Lambda$ such that either $x_I\neq
1$ or
$y_I\neq 1$, and write $x=x_Ix'$ and $y=y_Iy'$.  By Proposition \ref{lubs}
we may
write $x\vee y = (x_I\vee y_I)(x'\vee y')$, and hence
\[\phi(x\vee y)= (x_I\vee y_I)\phi(x'\vee y')\,.\] By induction on
$\ell(x)+\ell(y)$
we have that $\phi(x'\vee y')=\phi(x')\vee\phi(y')$. Thus it remains to
show that
\begin{equation}\label{laststep}
\phi(x)\vee \phi(y)= (x_I\vee y_I)\cdot(\phi(x')\vee \phi(y')).
\end{equation}

Note that in the direct product $\bigoplus_{J\in\Lambda} (G_J,P_J)$ every
element
$\xi$ may be written $\xi=\prod_{J\in\Lambda} \xi_J$. Combining this with
Proposition
\ref{lubs} we have:
\begin{equation}\label{laststep2}
\phi(x)\vee \phi(y)= (\phi(x)_I\vee\phi(y)_I)\cdot
\prod\limits_{J\in\Lambda ;J\neq I}(\phi(x)_J\vee\phi(y)_J)
\end{equation} Let Equation (\ref{laststep2}$'$) denote the equation similarly
obtained for the pair $x'$ and $y'$. We recall that $\phi(x)_J=\phi(x')_J$ and
$\phi(y)_J=\phi(y')_J$ for all $J\neq I$. Thus, Equation (\ref{laststep})
follows from
Equations (\ref{laststep2}) and (\ref{laststep2}$'$) via Lemma \ref{lemA}.
\Endproof

We can now extend  Theorem 6.7 of \cite{LR} to graph products of amenable
groups.

\begin{thm}\label{main-amen}
Any graph product of a family of quasi-lattice orders in
which the underlying groups are amenable is an amenable quasi-lattice order.
(i.e. the Toeplitz representation is faithful).
\end{thm}

\Proof Let $(G,P)$ denote a graph product of the family
$(G_I,P_I)_{I\in\Lambda}$ of quasi-lattice orders, and let $({\cal G},{\cal
P})$
denote their direct product.  By Proposition \ref{mainprop}, the map
$\phi:(G,P)\to({\cal G},{\cal P})$ induced by the identity on each factor
satisfies
the hypothesis of Proposition \ref{LacaRae}.  If each $G_I$ is amenable as
a group then
so is ${\cal G}$.  It then follows, by Proposition \ref{LacaRae}, that
$(G,P)$ is an
amenable quasi-lattice order.\Endproof

As an application of Theorem \ref{main-amen}, we obtain the following
characterisation
of faithfulness of representations of $C^*(G,P)$, and the consequent result
about
uniqueness of the $C^*$-algebra generated by a covariant isometric
representation of $P$. Monoid representations $V: P\to\Isom(\mathcal H)$ and $W:Q\to\Isom(\mathcal H)$ are said to \emph{$*$-commute} if every $V_x$ or $V_x^*$ for $x\in P$ commutes with every $W_y$ or $W_y^*$ for $y\in Q$, and to be \emph{orthogonal} to one another if $V_x^*W_y=0$ for all $x\in P$ and $y\in Q$.

\begin{thm}\label{unique-gr-pr}
Let $(G,P)$ be the graph product of a family
$(G_I,P_I)_{I\in \Lambda}$ of quasi-lattice ordered groups.

\begin{enumerate}
\item[(i)] If $\{ V_I: P_I \to \Isom(\mathcal H)\}_{I\in \Lambda}$ is a
family of covariant isometric representations
such that $V_I$ $*$-commutes with $V_J$ when $I$ and $J$ are adjacent in $\Gamma$ and
$ V_I$ is orthogonal to $V_J$ when $I$ and $J$ are not adjacent in $\Gamma$,
then there is a (unique) isometric covariant representation $V: P \to\Isom(\mathcal H)$ such that $V|_{P_I} = V_I$. All covariant representations
of $P$ arise this way.
\end{enumerate}

\noindent Suppose that each $G_I$, for $I\in\Lambda$, is an amenable group. 
Then we also have the following:  

\begin{enumerate}
\item[(ii)] The representation of
the
universal algebra $C^*(G,P)$ associated to $V$ is faithful if and only if
\begin{equation}\label{nontriv-cond}
\prod_{x\in F} (1 - V_x V_x^*) \neq 0 \qquad \text{ for every finite subset
} F\subset
\cup_{I\in \Lambda} (P_I \setminus \{1\})\,.
\end{equation}

\item[(iii)] If $\{V_I\}$ and $\{W_I\}$ are two families of covariant
isometric  representations as in part (i) satisfying
\eqref{nontriv-cond}, then the canonical map $V_I(x) \mapsto W_I(x)$
extends to an
isomorphism of $C^*(V_x: x\in P)$ onto $C^*(W_x: x\in P)$.
\end{enumerate}
\end{thm}

\Proof
The isometries satisfy the commuting relations which define $P$ as a graph product of monoids.
Therefore the maps $s\in P_I \mapsto V_I(s)$ extend to an isometric
representation $V$ of the monoid $P = \Gamma_{I\in \Lambda} P_I$.

We need to show that this representation is covariant, that is,
we need to show that
$V_{x\vee y}V_{x\vee y}^* = V_x V_x^* V_y V_y^*$
for all $x$ and $y$ in $P$.

We proceed by induction on $\ell(x) + \ell(y)$.
Choose $I\in \Delta(x) \cup \Delta(y)$
and write $x = x_I x'$ and $y = y_I y'$.
Since $x_I$ and $y_I$ cannot both be trivial,
this ensures that $\ell(x') + \ell(y') < \ell(x) + \ell(y)$.

Suppose first that $x\vee y < \infty$. By Proposition \ref{lubs}
and the induction hypothesis, we have
\begin{align}\label{zxxyyz}
V_{x\vee y}V_{x\vee y}^* = V_z V_{x'} V_{x'}^* V_{y'} V_{y'}^* V_z^*
\end{align}
where $z = x_I \vee y_I \in P_I$. Moreover, $x$ and $y$ are $I$-adjacent, and from Definition \ref{I-adj} (b) and (c) we have four
cases to consider:
\smallskip

\noindent{\sl Case 1:\ }
The vertex $I$ is adjacent to all the vertices of $x'$ and of $y'$. Then
$V_I$ $*$-commutes with $V_{x'}$ and with $V_{y'}$, so that \eqref{zxxyyz} becomes
\begin{align*}
 V_{x\vee y}V_{x\vee y}^* &= V_{x'} V_{x'}^* V_z V_z^* V_{y'} V_{y'}^* \\
&= V_{x'} V_{x'}^* V_{x_I} V_{x_I}^* V_{y_I} V_{y_I}^* V_{y'} V_{y'}^*
\quad \text{ by
covariance of } V_I\\
&= V_{x_I}V_{x'} V_{x'}^*  V_{x_I}^* V_{y_I}  V_{y'} V_{y'}^*V_{y_I}^*\\
&= V_x V_x^* V_y V_y^* \,.
\end{align*}

\noindent{\sl Case 2:\ }
 $z= x_I$ and $I$ is adjacent to all vertices of $y'$. Write
$z = y_I v$, with $ v\in P_I$; then
$V_{y'} V_{y'}^* V_z^* = V_v^* V_{y'} V_{y'}^*V_{y_I}^* = V_z^* V_{y} V_{y}^*$.
Thus \eqref{zxxyyz} becomes
\begin{align*}
V_{x\vee y}V_{x\vee y}^* & = V_z V_{x'} V_{x'}^* V_z^* V_{y} V_{y}^* \\
&= V_x V_x^* V_y V_y^* \,.
\end{align*}

\noindent{\sl Case 3:\ }
 $z= y_I$ and $I$ is adjacent to all vertices of $x'$.
(This is analogous to case 2.)

\noindent{\sl Case 4:\ } $z = x_I = y_I$.
Inserting $V_z^* V_z$ in the middle of the right-hand side of \eqref{zxxyyz} we get
$V_{x\vee y}V_{x\vee y}^* = V_x V_x^* V_y V_y^*$.

\medskip

Suppose now that $x\vee y = \infty$. Then by convention
$V_{x\vee y}V_{x\vee y}^* = 0$ and it suffices to show that
$V_x^* V_y =0$.
Clearly if $x_I \vee y_I = \infty$ then, by covariance of $V_I$, we have 
$V_{x_I}^* V_{y_I} =0$, so $V_{x}^* V_{y} =0$. Thus, we may suppose that 
$x_I\vee y_I\neq\infty$ and hence that
\[
V_{x_I}^* V_{y_I} = V_u V_v^* \qquad \text{where } x_I \vee y_I= x_I u = y_I v.
\]
Thus
\[
V_x^* V_y = V_{x'}^* V_u V_v^* V_{y'}.
\]

\noindent Let $A = a_1 a_2 \cdots a_k$ be a reduced expression for $y'$. Either 
condition (c) of Definition \ref{I-adj} holds, or $v\neq 1$ and $A$ has a 
syllable with vertex not adjacent to $I$. Suppose the latter, and let $a_i$ 
be the first syllable of $A$ for which $I(a_i)$ is not adjacent to $I$.
Note that $I(a_i)\neq I$ because $y'$ cannot have initial 
vertex $I$. Write 
$\alpha = a_1 a_2 \cdots a_{i-1} $ 
and $\beta =  a_{i+1} a_{i+2} \cdots a_{k}$.
We now have
\[
V_v^* V_{y'} = V_v^* V_\alpha V_{a_i} V_\beta = V_\alpha V_v^* V_{a_i}
V_\beta = 0
\]
using the fact that $V_v^* V_{a_i} = 0$ by orthogonality.

We may thus suppose that condition (c) of Definition \ref{I-adj} holds,
in which case $V_v^* V_{y'} = V_{y'} V_v^*$. By a similar argument we may 
suppose also that condition (b) holds and that $V_{x'}^* V_u  = V_u 
V_{x'}^*$. We have already assumed that $x_I\vee y_I\neq\infty$ (part (a) of 
Definition \ref{I-adj}). All these conditions together imply that $x$ and 
$y$ are $I$-adjacent (Definition \ref{I-adj}) and that
\[
V_x^* V_y = V_u V_{x'}^*  V_{y'} V_v^*.
\]
By Proposition \ref{lubs}(i) we now have $x'\vee y'=\infty$ (since $x\vee y=\infty$ and $x$ is $I$-adjacent to $y$), and applying the induction hypothesis to $x'$ and $y'$ completes the proof that $V$ is covariant.

\medskip

To prove (ii) observe first that from Theorem \ref{main-amen} and
\cite[Theorem 3.7]{LR}, it follows that a representation is faithful if and
only if
\[
\prod_{x\in F} (1 - V_x V_x^*) \neq 0 \qquad \text{ for every finite subset
} F\subset P.
\]
It is equivalent to consider only products of the form stated in
\eqref{nontriv-cond} because replacing each $x \in F$ by one of its initial
syllables
has the effect of replacing each factor $(1- V_x V_x^*)$ by a smaller one. Part
(iii) follows from Theorem \ref{main-amen} and Corollaries 3.8 and 3.9 of \cite{LR}.
\Endproof

\begin{remark} In some cases  condition \eqref{nontriv-cond} is
automatically satisfied
by all covariant representations, in which case the Toeplitz $C^*$-algebra
is simple,
and purely infinite by \cite[Theorem 5.4]{purelinf}. The best known example
of this is
$\mathcal O_\infty$ \cite{cun}. See Corollaries 5.2 and 5.3 of \cite{LR}
and  Theorem
2.4 of \cite{dav-pop} for more examples involving free products.
\end{remark}

\section{The $C^*$-algebra of a right-angled Artin semigroup of isometries.}

Let $\Lambda$ be a set (usually taken to be finite, although we shall not
make this
restriction here). A matrix $\Mat = (m_{s,t})_{s,t\in \Lambda}$ is a {\em
Coxeter
matrix} if
 $m_{s,t}= m_{t,s} \in \{2,3,\ldots, \infty\}$ for $s\neq t$ and
$m_{s,s} =1$.
 Denote by
$\langle s t \rangle^m$ the word  $sts \cdots$,  beginning with $s$ and
having length
$m$, in which the letters $s$ and $t$ alternate.

The {\em Artin group} $A_\Mat$ associated to $\Mat$ is the group with
presentation

\[
\langle \ \Lambda\ | \ \langle s t\rangle^{m_{s,t}} = \langle t
s\rangle^{m_{s,t}}\text{ for each }  s,t \in \Lambda\ \rangle,
\]

\noindent in which a relation of the form $\langle s t\rangle^\infty =
\langle t s\rangle^\infty$ is to be interpreted as vacuous.  The {\em Artin
monoid}
$A_\Mat^+$ is defined via the same presentation, taken in the category of
monoids
(semigroups with unit), see \cite{BS}.  We may view $A_M$ as a partially
ordered group
with positive cone $P_M$ generated by $\Lambda$.  The cone $P_M$ is in
general a
quotient of
$A_M^+$ via the obvious map, although in many cases of interest this map is
known to
be injective (see \cite{Charney} for the most recent results).

Adding the relations $s^2 =1$ for $s\in \Lambda$ to the above ones yields a
somewhat
unusual presentation of the more familiar Coxeter group associated with
$\Mat$, which
is usually presented via the relations
$(st)^{m_{s,t}} = 1$.

\begin{defn} A Coxeter group and its associated Artin group $A_\Mat$ are
said to be
{\em right-angled} if every nondiagonal entry of the Coxeter matrix $\Mat$
is either
$2$ or $\infty$; see for example
\cite{Chiswell}.  (By abuse we also refer to the matrix as right-angled.)
\end{defn}

The terminology is motivated by noting that the right-angled Coxeter groups
are those
linear reflection groups whose reflecting hyperplanes are mutually
orthogonal or
parallel.

Every right-angled Coxeter matrix $\Mat$ over $\Lambda$ determines a graph
$\Gamma$
with vertex set $\Lambda$ having an edge joining $s$ and
$t$ when $m_{s,t} = 2$.  The only relations in the presentation of
$A_\Mat$ say that two generators commute if they are joined by an edge,
hence $A_\Mat$
is precisely the graph product $\Gamma_{I\in \Lambda}
\mathbb Z$ of copies of $\mathbb Z$.  In this connection, right-angled
Artin groups
are also referred to as \emph{graph groups}, see \cite{HM}.

By virtue of Theorem \ref{NormalForm} the right-angled Artin monoid
$A_\Mat^+$ may be identified with the positive cone of the corresponding
Artin group
$A_\Mat$, and  $(A_\Mat, A_\Mat^+)$ is a quasi-lattice order by Theorem
\ref{graph-prods-qlo}.  Applying  Theorem \ref{main-amen} we see that this
quasi-lattice order is amenable, and hence the Toeplitz representation of
$C^*(A_\Mat,
A_\Mat^+)$ is faithful. As with Coburn's and Cuntz's theorems, it is more
appealing to
formulate the result in terms of the generators themselves; indeed, notice that
assertion (iii) below does not contain any explicit reference to
quasi-lattice orders
or Artin groups.

\begin{thm} \label{ra-artin}
Let $\Gamma$ be a graph with set of vertices $\Lambda$
and suppose $\{V_s: s\in \Lambda\}$ is a collection of isometries on
Hilbert space
such that for every pair of distinct vertices $s$ and $t$ one has
\[
V_s V_t = V_t V_s \quad \text{ and } \quad V_s^* V_t = V_t V_s^* \quad
\text{ if $s$ and $t$ are adjacent in } \Gamma \text{, and}
\]
\[
V_s^* V_t = 0 \quad \text{ if $s$ and $t$ are not adjacent in } \Gamma.
\]
 Let $A_\Mat = \Gamma_{s \in \Lambda} \mathbb Z\< s\>$ be the right-angled
Artin group
associated to $\Gamma$. Then
\begin{enumerate}
\item[(i)]  the maps $s\to V_s$, for $s\in\Lambda$,
extend to a covariant isometric representation $V$ of the right-angled Artin
semigroup
$A_\Mat^+$,

\item[(ii)] the corresponding representation of $C^*(A_\Mat, A_\Mat^+)$ is
faithful if and only if
\begin{equation}\label{nondeg-ra-artin}
\prod_{s\in F} (I - V_s V_s^*) \neq 0 \qquad \text{ for every finite }F
\subset \Lambda, \quad \text{ and}
\end{equation}
\item[(iii)] if $\{W_s: s\in \Lambda\}$ is another collection of isometries
satisfying
the same relations  and condition  \eqref{nondeg-ra-artin}, then the map $V_s
\mapsto W_s$ extends to a $C^*$-algebra isomorphism of $C^*(V_s: s\in
\Lambda)$ to
$C^*(W_s: s\in \Lambda)$.
\end{enumerate}
\end{thm}

\Proof Since the generators $V_s$ satisfy the stated relations, the
collection of
isometric representations $\{n\in \N  \mapsto  V_s^n\}_{s\in \Lambda}$
satisfies the hypothesis of Theorem \ref{unique-gr-pr} (i) and so extends to an
isometric representation of the semigroup $A_\Mat^+$, giving (i). The
amenability
hypothesis is satisfied because in this case each factor in the graph product is isomorphic to
$\Z$, so (ii)
and (iii) follow directly from Theorem \ref{unique-gr-pr}.
Notice that the necessary and sufficient condition for faithfulness in (ii)
is equivalent to the one obtained in Theorem \ref{unique-gr-pr} because every
projection of the form \eqref{nontriv-cond} majorates a projection of the form
\[
\prod_{s\in F} (I - V_s V_s^*)\qquad \text{ for some finite }F \subset
\Lambda,
\]
to see this it suffices to replace each syllable $x=s^n\in \N\< s\>$ by the
corresponding generator $s$.
\Endproof

\begin{remark} If the set $\Lambda$ of generators is finite, the projection
$\prod_{s\in \Lambda} (I - T_s T_s^*)$ belongs to the Toeplitz algebra
$\mathcal T (A_M, A_M^+)$. Since every nontrivial element in $A_M^+$ is
bounded below
by a generator, it follows from \cite[Proposition 6.3]{Nica} that the compact
operators on $\ell^2(P)$ are contained in $\mathcal T (A_M, A_M^+)$ as the
ideal
generated by this projection. It is easy to see that
\eqref{nondeg-ra-artin} holds if
and only if it holds for $F = \Lambda$, so this ideal is minimal.
\end{remark}

\section{Non-amenability of lattice ordered groups}\label{nonamen}

Next we concentrate on partially ordered groups in which least comon upper
bounds
always exist.

\begin{defn} A partially ordered group $(G,P)$ is \emph{lattice ordered} if
every pair
of elements has a least common upper bound. By left invariance, an equivalent
condition is that every element $x\in G$ have a least upper bound in $P$.
\end{defn}

Lattice orders are special cases of quasi-lattice orders; in fact we have the
following characterisation, cf.~\cite{Nica}.

\begin{lemma} The following are equivalent for a partially ordered group
$(G,P)$:

\begin{enumerate}
\item[(i)] $(G,P)$ is lattice ordered.

\item[(ii)] $(G,P)$ is quasi-lattice ordered and $G = PP\inv$.

\item[(iii)] $(G,P)$ is quasi-lattice ordered, $P$ generates $G$
             (in which case we say $(G,P)$ is connected),
              and $aP\cap bP\neq \emptyset$ for all $a,b\in P$.
\end{enumerate}
  \end{lemma}

\Proof
 (i)$\implies$(ii): Suppose $(G,P)$ is lattice ordered and let $x\in G$;
then  $x\leq
x\vee 1$ so that both $a := x\vee 1$ and $b := x\inv a$ are in $P$. Clearly
$x =
ab\inv$.

(ii)$\implies$(iii): If $G=PP\inv$ then obviously $P$ generates $G$. Given
a pair $a,
b \in P$, write $a\inv b$ as $x y\inv$ with $x,y \in P$. Then $ax = by \in
aP \cap bP$.

(iii)$\implies$(i): Suppose (iii) holds. Take any $x\in G$ and (using the connectedness) write $x=a_1a_2...a_k$ with each
$a_i\in P$ or $P\inv$. Since any element $a\inv b$ for $a,b\in P$ may be rewritten $cd\inv$ for $c,d\in P$ (by finding $ac = bd\in aP\cap bP$), any such expression for $x$ may ultimately be simplified to the form $x=uv\inv$ with $u,v\in P$. Thus $x$ has an upper bound $u\in P$ and, 
since $(G,P)$ is a quasi-lattice order, it must have a least
upper bound in $P$. So (i) holds.\Endproof

For a quasi-lattice order $(G,P)$, we know that if $G$ is amenable then
$(G,P)$ is
amenable, by \cite[\S 4.5]{Nica}, see also \cite[Lemma 6.5]{LR}. It turns
out that for
lattice ordered groups the converse is also true; the proof follows the
argument
outlined in Remark 2 of \cite[\S 5.1]{Nica}.

\medskip
\begin{prop} \label{latt-amen}
 If $(G,P)$ is lattice ordered and amenable (in the sense of Definition
\ref{defqlo}) then $G$ is an amenable group.
\end{prop}

\Proof Suppose $(G,P)$ is an amenable lattice order, and denote the left
regular
representation of $P$ on $\ell^2(P)$  by $W$. As observed by Nica, the map
$x\in P \mapsto 1 \in \mathbb C$ is a (one-dimensional) covariant isometric
representation, so the map $W_x \mapsto 1$ extends to give a one-dimensional
representation of
$\mathcal T(G,P)$, and an easy argument shows that
\begin{equation}\label{nicaell1}
\| \sum_{s\in P} \lambda_s W_s \| = \sum_{s\in P} \lambda_s
\end{equation} for every finitely (or countably) supported  nonnegative
function
$\lambda$ in $\ell^1(P)$.

Notice  that the sum $\sum_s \lambda_s W_s$ is the operator of
left-convolution by
$\lambda \in \ell^1(P)$ on $\ell^2(P)$, and that \eqref{nicaell1} implies
that $\|
\sum_{s\in P}
\lambda_s W_s \| = 1$ for every probability density $\lambda$. Since the
support of
$\lambda$ can be chosen to contain an arbitrary finite subset of $P$
including the
identity, we have that condition (e) of \cite[Theorem 1]{Day} holds, with
$p = 2$, $U = 1 \in P$, and $\phi $ a probability density whose support
contains $\xi
\cup \{1\}$. By \cite[Theorem 1]{Day}, the semigroup  $P$ has a
left-invariant mean.
(Day's Theorem is about right-amenability and right-convolutions, but there
is no
difficulty in transforming it into a theorem for left-amenability and
left-convolutions.)

Finally, the group $G = P P\inv$ is amenable, by Corollary 3.6 of
\cite{WiWi}.\Endproof

The Artin group $A_\Mat$ is said to be of {\em finite type} if $\Lambda$ is
finite and
the Coxeter group associated to the matrix $\Mat$ is finite, \cite{BS}.
 By Proposition 5.5 and Theorem 5.6 of \cite{BS}, if $A_\Mat$ is of finite
type, then
the Artin semigroup $A_\Mat^+$ embeds as a subsemigroup of $A_\Mat$, and
the pair
$(A_\Mat, A_\Mat^+)$ is a lattice ordered group (see also \cite{Del}). We
wish to
apply Proposition \ref{latt-amen} to characterize amenability of these
lattice orders;
the first observation is that most Artin groups are not amenable.

\begin{prop}\label{artin-amen}  Let $M$ be a Coxeter matrix over a finite set
$\Lambda$. The Artin group $A_M$ is amenable (as a group) if and only if it
is the
free abelian group on $\Lambda$, that is $m_{s,t}=2$ for all
$s,t\in\Lambda$, $s\neq
t$.
\end{prop}

\Proof In \cite{CP} it is shown that the set of elements
$\Cal Q = \{ Q_s=s^2 :s\in\Lambda\}$
generates a subgroup $H_M$ of $A_M$ with presentation
\begin{equation}\label{HMpres}
\langle\Cal Q\mid Q_sQ_t=Q_tQ_s\text{ if } m_{s,t}=2\rangle\,.
\end{equation}
That is, $H_M$ is a right-angled Artin group. In particular, $H_M$ contains 
at least one free subgroup of rank 2, and hence is not amenable, unless of 
course $m_{s,t}=2$ for all $s\neq t\in\Lambda$. In the latter case $A_M$ is 
free abelian and therefore amenable.\Endproof

As a consequence of Proposition \ref{artin-amen}, the analog of Theorem
\ref{unique-gr-pr} (ii) fails for nonabelian Artin groups of finite type,
instead we have the following non-uniqueness result.

\medskip
\begin{thm}\label{nonunique} Let $A_\Mat$ be a nonabelian Artin group of
finite type.
Suppose $\{V_s: s\in \Lambda\}$ is a collection of isometries satisfying
the Artin
relations relative to $M$:
\begin{equation}\label{rels-ftype}
\langle V_s V_t\rangle^{m_{s,t}} = \langle V_t V_s\rangle^{m_{s,t}} \qquad
s,t\in
\Lambda \end{equation} Then the map $s \mapsto V_s$ extends to an isometric
representation, denoted also by $V$, of the Artin semigroup $A_M^+$. The
representation $V$ is covariant provided that
\begin{equation}\label{cov-ftype} V_s V_s^* V_t V_t^* = V_{s\vee t}
V_{s\vee t}^*
\qquad s,t\in\Lambda\,.
\end{equation}

The Toeplitz representation $T$ on $\ell^2(A_M^+)$ satisfies
\eqref{rels-ftype} and
\eqref{cov-ftype} and, moreover, the projection
$\prod_{\Lambda} (1 - T_s T_s^*) $ does not vanish. However, the Toeplitz
representation of $C^*(A_\Mat,A_\Mat^+)$ is not faithful; in particular, the
$C^*$-algebra generated by a collection $\{V_s: s\in \Lambda\}$ as above is not
canonically unique, even if we assume
$ \prod_{\Lambda} (1 - V_s V_s^*) \neq 0 $.
\end{thm}

\Proof The map can be extended to an isometric representation because the given
isometries satisfy the relations
\eqref{rels-ftype}, which constitute a presentation of $A_M^+$.

Suppose that (\ref{cov-ftype}) holds.  We need to show that the covariance
condition
holds for every pair $x,y \in A_M^+$.  We use the length homomorphism
$l:A_M^+\to \N$
such that $l(s)=1$ for each generator $s\in\Lambda$.  Choose $x,y$ such
that the
covariance condition is not satisfied and such that $l(x\vee y)$ is
minimised. Amongst
the possible pairs, choose one so that $l(x)$ is the smallest possible.
\smallskip

\noindent {\sl Case 1:} Suppose $l(x)=1$, so that $x$ is actually a
generator $s$.
Write $y=tz$ for $t$ a generator, so that $V_y V_y^* = V_t(V_t^*V_t)V_zV_y^* =
V_tV_t^* V_y V_y^*$. Then, by (\ref{cov-ftype}),
\[ V_s V_s^* V_y  V_y^*= V_s V_s^* V_tV_t^* V_y V_y^* = V_{s\vee t}
V_{s\vee t}^* V_y
V_y^*.
\] Now $s\vee t$ and $y$ have a common left factor $t$; writing
$s\vee t = tu $ and $y = tz$ we have
\[ V_{s\vee t} V_{s\vee t}^* V_y V_y^* = V_t V_uV_u^* V_z V_z^* V_t^*.
\] Since  $s\vee y = s \vee t \vee y = t (u\vee z)$, we have that
$l(u\vee z)<l(s\vee y)$ and  covariance follows by applying the induction
hypothesis
to $V_u V_u^* V_z V_z^*$ in this expression.
\smallskip

\noindent {\sl Case 2:} Suppose now $l(x)>1$, and write $x=su$ for $s$ a
generator.
Put $s\vee y=sz$, so that $x\vee y = x\vee s\vee y = s(u\vee z)$. Since
$l(s\vee
y)\leq l(x\vee y)$ and $l(s)<l(x)$, the induction hypothesis implies that
\[ V_s V_s^* V_y  V_y^* = V_s V_z V_z^* V_s^*\,.
\] Applying $V_s^*$ on the left and using the fact that $V_s^*V_s=1$ this
equation
becomes $V_s^* V_y V_y^*=V_z V_z^* V_s^*$. We then have:
\[ V_x V_x^* V_y  V_y^*= V_s V_u V_u^* V_s^* V_y  V_y^* = V_s V_u V_u^* V_z
V_z^*
V_s^*\,.
\] Since $x\vee y=s(u\vee z)$ we have $l(u\vee z) < l(x\vee y)$, and 
covariance follows in this case, as before, by applying the induction 
hypothesis to $V_u V_u^* V_z V_z^*$.
\smallskip

It is clear that the Toeplitz representation satisfies the stated
relations, and the
projection $\prod_{\Lambda} (1 - T_s T_s^*)$ does not vanish at
$\delta_1 \in \ell^2(A_M^+)$.
The quasi-lattice ordered group $(A_M, A_M^+)$ is not amenable by Proposition
\ref{latt-amen} and Proposition \ref{artin-amen} and the last assertion
follows, cf.
\cite[Corollary 3.9]{LR}. \Endproof

\begin{remark} Since Theorem \ref{nonunique} shows that the Toeplitz
algebra $\mathcal
T(A_M,A_M^+)$ of a nonamenable finite type Artin group 
$A_M$ is not universal for covariant isometric representations, it is
generally hard
to decide whether a given collection of isometries satisfying the Artin
relations and
the covariance condition actually generates a representation of $\mathcal
T(A_M,A_M^+)$. In any case, it follows from Theorem 6.7 and Corollary 6.8
of \cite{elq}
that
 a given a representation of $\mathcal T(A_M,A_M^+)$ is faithful if and
only if the
generating family of isometries is proper, in the sense that
$\prod_{\Lambda} (1 - V_s V_s^*) \neq 0 $.
\end{remark} Our results about Toeplitz algebras cover the Artin groups
$(A_M,A^+_M)$ that are presently known to be quasi-lattice ordered, namely
the finite
type Artin groups from \cite{BS,Del} and the right-angled Artin groups
dealt with by
Theorem \ref{graph-prods-qlo}. It would be interesting to formulate and decide
questions of amenability and uniqueness in the remaining cases. It is known
that, for
$\Lambda$ finite, the monoid $A_M^+$ always has a quasi-lattice structure (see
\cite{BS, Del}), but even when it is known that $A_M^+$ embeds canonically in
$A_M$ (which remains an open question in general; see \cite{Charney}) this
is not
enough to show that $(A_M,A_M^+)$ is quasi-lattice ordered, which is
essential for our
techniques.

\vskip1cm

John Crisp
\medskip

Laboratoire de Topologie

Universit\'e de Bourgogne

UMR 5584 du CNRS, B.P. 47 870

21078 Dijon Cedex

FRANCE

\smallskip

\texttt{crisp@topolog.u-bourgogne.fr}

\vskip1.5cm

Marcelo Laca
\medskip

Department of Mathematics

The University of Newcastle

NSW 2308

AUSTRALIA

\smallskip

\texttt{marcelo@math.newcastle.edu.au}

\end{document}